\documentclass{article}
\usepackage{amsmath, amsthm}
\usepackage{amssymb}
\newtheoremstyle{theorem}
  {10pt}          
  {10pt}  
  {\sl}  
  {\parindent}     
  {\bf}  
  {. }    
  { }    
  {}     
\newtheorem{thm}{Theorem}[section]
\newtheorem{cor}[thm]{Corollary}
\newtheorem{lem}[thm]{Lemma}

\newtheorem{defn}{Definition}[section]
\newtheorem{exam}[defn]{Example}
\newtheorem{rem}{Remark}[section]
\begin{document}
\title{\Large\bf The \emph{q}--fractional analogue for Gronwall--type
inequality}

\author{
Thabet Abdeljawad$^{a,b}$,
Jehad O. Alzabut$^{b}$\\
\small{$^a$Department of Mathematics and Computer Science, \c{C}ankaya University}\\
\small{\"{O}gretmenler Cad. 14, 06530, Balgat--Ankara, Turkey}\\
\small{thabet@cankaya.edu.tr}\\
\small{$^{b}$Department of Mathematics and Physical Sciences, Prince Sultan University}\\
\small{P.O.Box 66833, Riyadh 11586, Saudi Arabia}\\
\small{jalzabut@psu.edu.sa}\\
}
\date{}
\maketitle {\small{ \noindent {\bf Abstract.} In this article, we
utilize \emph{q}--fractional Caputo initial value problems of order
$0<\alpha\leq 1$ to derive a \emph{q}--analogue for Gronwall--type
inequality. Some particular cases are derived where
\emph{q}--Mittag--Leffler functions and \emph{q}--exponential type
functions are used. An example is given to illustrate the
validity of the derived inequality.}

\vskip.1in

\section{Introduction} \label{s:1}
The fractional differential equations have conspicuously received
considerable attention in the last two decades. Many researchers
have investigated these equations due to their  significant
applications in various fields of science and engineering such as
in viscoelasticity, capacitor theory, electrical circuits,
electro--analytical chemistry, neurology, diffusion, control
theory and statistics; see for instance the monographs
\cite{Samko,Pody,Kilbas}. The study of \emph{q--}difference
equations, on the other hand, has gained intensive interest in the
last years. It has been shown that these types of equations have
numerous applications in diverse fields and thus have  evolved
into multidisciplinary subjects \cite{F1,F2,F3,F4,F5,F6,F7}. For
more details on \emph{q--}calculus, we refer the reader to the
references \cite{F8,F9}. The corresponding fractional difference
equations, however, have been comparably less considered. Indeed,
the notions of fractional calculus and \emph{q--}calculus are
tracked back to the works of Euler  and Jackson \cite{J1},
respectively. However, the idea of fractional difference equations
is considered to be very recent; we suggest the new papers \cite
{T1,T2,T3,T4,T5,T6',T6,T7',T7,T8,T9,T10,T11,JC1,JC2} whose authors
have taken the lead to promote the theory of fractional difference
equations.

The  \emph{q}--fractional difference equations which serve as a
bridge between fractional difference equations and
\emph{q}--difference equations have become a main object of
research in the last years. Recently, there have appeared many
papers which study the qualitative properties of solutions for
\emph{q}--fractional differential equations
\cite{AT1,AT2,AT3,AT4,AT5} whereas few results exist for \emph{q}--fractional
difference equations \cite{AT6,AT7,AT8}.  The integral
inequalities which are considered as an effective tools for
studying solutions properties have been also under consideration.
In particular, we are interested with Gronwall's inequality which has been a main target for
many researchers. There are several versions for Gronwall's
inequality in the literature; we list here those results which
concern with fractional order equations
\cite{d2,d3,d4,d5,d1}. To the best of authors' observation, however, the
\emph{q}--fractional analogue for Gronwall--type inequality has
not been investigated yet.

A primary purpose of this paper is to utilize the
\emph{q}--fractional Caputo initial value problems of order
$0<\alpha\leq 1$ to derive a \emph{q}--analogue for Gronwall--type
inequality. Some particular cases are derived where
\emph{q}--Mittag--Leffler functions and \emph{q}--exponential type
functions are used. An example is given to illustrate the validity
of the derived inequality.

\section {Preliminary assertions}\label{s:3}
Before stating and proving our main results, we introduce some definitions and
notations that will be used throughout the paper. For $0<q<1$, we define  the time scale  $\mathbb{T}_q$ as follows
$$\mathbb{T}_q=\{q^n:n \in \mathbb{Z}\}\cup \{0\}, $$ where $\mathbb{Z}$ is the set of
integers. In general, if $\alpha$ is a nonnegative real number
then we define the time scale
$$\mathbb{T}_q^\alpha=\{q^{n+\alpha}:n \in \mathbb{Z}\}\cup \{0\}$$ and thus we may
write
$\mathbb{T}_q^0=\mathbb{T}_q.$ For a function $f:\mathbb{T}_q\rightarrow
\mathbb{R}$, the nabla
$q-$derivative of $f$ is given by
\begin{equation} \label{qd}
\nabla_q f(t)=\frac{f(t)-f(qt)}{(1-q)t},~~t \in
\mathbb{T}_q-\{0\}.
\end{equation}
The nabla \emph{q}--integral of $f$ is given by
\begin{equation} \label{qi}
\int_0^t f(s)\nabla_q s=(1-q)t\sum_{i=0}^\infty q^if(tq^i)
\end{equation}
and
\begin{equation} \label{qiyyyy}
\int_a^t f(s)\nabla_q s=\int_0^t f(s)\nabla_q s - \int_0^a
f(s)\nabla_q s,\;\;\mbox{for}\;\;0\leq a \in T_q.
\end{equation}
The \emph{q}--factorial function for $n\in \mathbb{N}$ is defined
by
\begin{equation} \label{qfact}
(t-s)_q^n=\prod_{i=0}^{n-1}(t-q^is).
\end{equation}
In case $\alpha$ is a non positive integer, the
\emph{q}--factorial function is defined by
\begin{equation} \label{qfactg}
(t-s)_q^\alpha=t^\alpha\prod_{i=0}^\infty \frac{1- \frac{s}{t}
q^i} {1- \frac{s}{t} q^{i+\alpha}}.
\end{equation}
In the following lemma,  we present some properties of
\emph{q}--factorial functions.
\begin{lem}\emph{\cite{AT4}}\label{qproperties} For $\alpha,\gamma,\beta \in
\mathbb{R}$, we have
\begin{itemize}
    \item [\emph{I.}] $(t-s)_q^{\beta+\gamma}=(t-s)_q^\beta (t-q^\beta
    s)_q^\gamma$.
    \item [\emph{II.}] $(at-as)_q^\beta=a^\beta (t-s)_q^\beta$.
    \item [\emph{III.}] The nabla q--derivative of the q--factorial function with
respect to $t$ is
$$\nabla_q (t-s)_q^\alpha
=\frac{1-q^\alpha}{1-q}(t-s)_q^{\alpha-1}.$$
    \item [\emph{IV.}] The nabla q--derivative of the q--factorial function with respect
to $s$ is
$$\nabla_q (t-s)_q^\alpha
=-\frac{1-q^\alpha}{1-q}(t-qs)_q^{\alpha-1}.$$
\end{itemize}
\end{lem}
For a function $f:\mathbb{T}^\alpha_q\rightarrow \mathbb{R}$, the
left \emph{q}--fractional integral $_{q}\nabla_a^{-\alpha}$ of
order $\alpha\neq 0,-1,-2,\ldots$ and starting at $0 <a \in
\mathbb{T}_q$ is defined by
\begin{equation} \label{ag}
 _{q}\nabla^{-\alpha}_{a} f(t)=\frac{1}{\Gamma_q(\alpha)}\int_a^t
 (t-qs)_q^{\alpha-1}f(s)\nabla_qs,
\end{equation}
where
\begin{equation} \label{qid}
\Gamma_q(\alpha+1)=\frac{1-q^\alpha}{1-q}\Gamma_q(\alpha),~~\Gamma_q(1)=1,~\alpha
>0.
\end{equation}
One should note that the left \emph{q}--fractional integral
$_{q}\nabla^{-\alpha}_{a}$ maps functions defined on
$\mathbb{T}_q$ to functions defined on $\mathbb{T}_q$.
\begin{defn} \emph{\cite{T1}}
If $0<\alpha \notin \mathbb{N}$ . Then the
 Caputo left q--fractional derivative of order $\alpha$  of a function
$f$ is  defined by
\begin{equation} \label{qrd}
_{q}C_a^\alpha f(t): = _{q}\nabla_a ^{-(n-\alpha)}\;\;\nabla_q
^nf(t)=\frac{1}{\Gamma(n-\alpha)} \int_a^t(t-qs)_q^{n-\alpha-1}
\nabla_q^nf(s)\nabla_q s,
\end{equation}
where $n=[\alpha]+1$. In case $\alpha \in \mathbb{N}$,  we may
write $_{q}C_a^\alpha f(t):= \nabla_q^n f(t)$.
\end{defn}
\begin{lem} \label{qtrans}\emph{\cite{T1}}
Assume that $\alpha>0$ and $f$ is defined in a suitable domain.
Then
\begin{equation}\label{qtrans1}
_{q}\nabla^{-\alpha}_{a}~ _{q}C_a^\alpha
f(t)=f(t)-\sum_{k=0}^{n-1}\frac{(t-a)_q^{k}}{\Gamma_q(k+1)}\nabla_q^kf(a)
\end{equation}
and if $0<\alpha\leq1$ then
\begin{equation}\label{qtrans3}
_{q}\nabla^{-\alpha}_{a}~ _{q}C_a^\alpha f(t)= f(t)-f(a).
 \end{equation}
\end{lem}
For solving linear \emph{q}--fractional equations, the following
identity  is essential
\begin{equation}\label{qpower}
_{q}\nabla^{-\alpha}_{a} (x-a)_q^\mu
=\frac{\Gamma_q(\mu+1)}{\Gamma_q(\alpha+\mu+1)}(x-a)_q^{\mu+\alpha},~~0<a<x<b,
\end{equation}
where $\alpha \in \mathbb{R}^+$ and $\mu \in (-1,\infty)$.  See
for instance the recent paper \cite{TD2011} for
more information.\\
\\
The \emph{q}--analogue of Mittag--Leffler function with double
index $(\alpha,\beta)$ is first introduced in \cite{TD2011}.
Indeed, it was defined as follows:
\begin{defn}\label{Mitt}\emph{\cite{T1}}
For $z ,z_0\in \mathbf{C}$ and $\mathfrak{R}(\alpha)> 0$, the
q--Mittag--Leffler function is defined by
\begin{equation} \label{Mit}
_{q}E_{\alpha,\beta}(\lambda,z-z_0)=\sum_{k=0}^\infty \lambda^k
\frac{(z-z_0)_q^{\alpha k}}{\Gamma_q(\alpha k+\beta)}.
\end{equation}
In case $\beta=1$, we may use
$~_{q}E_{\alpha}(\lambda,z-z_0):=~_{q}E_{\alpha,1}(\lambda,z-z_0)$.
\end{defn}
The following example clarifies how  \emph{q}--Mittag--Leffler
functions can be used to express the solutions of Caputo
\emph{q}--fractional linear initial value problems.
\begin{exam}\emph{\cite{T1}} \label{qlinear}
Let $0<\alpha\leq 1$ and consider the left Caputo
q--fractional difference equation
\begin{equation} \label{lfractional}
~_{q}C^ \alpha _a y(t)= \lambda y(t)+f(t),~~y(a)=a_0,~t\in T_q.
\end{equation}
Applying $_{q}\nabla^{-\alpha} _a$ to equation (\ref{lfractional})
and using (\ref{qtrans3}), we end up with
\begin{equation}\label{rep}
y(t)= a_0+ \lambda~ _{q}\nabla^{-\alpha} _a
y(t)+~_{q}\nabla^{-\alpha} _a f(t).
\end{equation}
To obtain an explicit form for the solution, we apply the method
of successive approximation. Set $y_0(t)=a_0$ and
$$y_m(t)=a_0+\lambda~ _{q}\nabla^{-\alpha} _a y_{m-1}(t)+
_{q}\nabla^{-\alpha} _a f(t), m=1,2,3,\ldots. $$ For $m=1$, we
have by the power formula (\ref{qpower})
 $$y_1(t)=a_0[1+\frac{\lambda ( t-a)_q^{(\alpha)}}{\Gamma_q(\alpha+1)}]+~
_{q}\nabla^{-\alpha} _a f(t).$$ For $m=2$, we also see that
\begin{eqnarray*}y_2(t)&=& a_0+ \lambda a_0~_{q}\nabla^{-\alpha} _a \Big[1+
\frac{(t-a)_q^\alpha}{\Gamma_q(\alpha+1)}\Big]+~_{q}\nabla^{-\alpha}
_a f(t) +\lambda~ _{q}\nabla_a^{-2\alpha} f(t)\\
&=&a_0 \Big[1+\frac{\lambda
(t-a)_q^\alpha}{\Gamma_q(\alpha+1)}+\frac{\lambda^2
(t-a)_q^{2\alpha}}{\Gamma_q(2\alpha+1)}\Big]+~_{q}\nabla^{-\alpha}
_a f(t) +\lambda~ _{q}\nabla_a^{-2\alpha} f(t).
\end{eqnarray*}
If we proceed inductively and let $m\rightarrow\infty$, we obtain
the solution
\begin{eqnarray*}y(t)&=&a_0 \Big[1+\sum_{k=1}^\infty\frac{\lambda^k (t-a)_q^{k\alpha}}
{\Gamma_q(k\alpha+1)}\Big]+\int_a^t
\Big[\sum_{k=1}^{\infty}\frac{\lambda^{k-1}}{\Gamma_q(\alpha
k)}(t-qs)_q^{\alpha k-1}\Big]f(s)\nabla_qs\\
&=& a_0\Big[1+\sum_{k=1}^\infty\frac{\lambda^k (t-a)_q^{k\alpha}}
{\Gamma_q(k\alpha+1)}\Big]+\int_a^t
\Big[\sum_{k=0}^{\infty}\frac{\lambda^{k}}{\Gamma_q(\alpha
k+\alpha)}(t-qs)_q^{\alpha k+(\alpha-1)}\Big]f(s)\nabla_qs\\
&=& a_0\Big[1+\sum_{k=1}^\infty\frac{\lambda^k (t-a)_q^{k\alpha}}
{\Gamma_q(k\alpha+1)}\Big]+\int_a^t
(t-qs)_q^{(\alpha-1)}\Big[\sum_{k=0}^{\infty}\frac{\lambda^{k}}{\Gamma_q(\alpha
k+\alpha)}(t-q^\alpha s)_q^{(\alpha k)}\Big]f(s)\nabla_qs.
\end{eqnarray*}
That is,
$$y(t)=a_0 ~_{q}E_\alpha (\lambda,t-a)+ \int_a^t
(t-qs)_q^{\alpha-1}~_{q}E_{\alpha,\alpha} (\lambda,t-q^\alpha s)
f(s)\nabla_q s.$$
\end{exam}
\begin{rem}If instead we use the  modified q--Mittag--Leffler function
$$~_{q}e_{\alpha, \beta}(\lambda,z-z_0)=\sum_{k=0}^\infty \lambda^k
\frac{(z-z_0)_q^{\alpha k+(\beta-1)}}{\Gamma_q(\alpha k+\beta)}$$
then, the solution representation (\ref{rep}) becomes
$$y(t)=a_0 ~_{q}e_\alpha (\lambda,t-a)+ \int_a^t
~_{q}e_{\alpha,\alpha} (\lambda,t-q s) f(s)\nabla_q s.$$
\end{rem}
\begin{rem}
If we set $\alpha=1$, $\lambda=1$, $a=0$ and $f(t)=0$, we reach to
the q--exponential formula  $e_q(t)= \sum_{k=0}^\infty
\frac{t^k}{\Gamma_q(k+1)}$ on the time scale $\mathbb{T}_q$, where
$\Gamma_q(k+1)=[k]_q!=[1]_q[2]_q...[k]_q$ with
$[r]_q=\frac{1-q^r}{1-q}$. It is known that $e_q(t)=E_q((1-q)t)$,
where $E_q(t)$ is a special case of the basic  hypergeometric
series, given by
$$E_q(t)=~_{1}\phi_0(0;q,t)=\prod_{n=0}^\infty
(1-q^nt)^{-1}=\sum_{n=0}^\infty \frac{t^n}{(q)_n},$$ where
$(q)_n=(1-q)(1-q^2)...(1-q^n)$ is the q--Pochhammer symbol.
\end{rem}
\section{The Main Results}
Throughout the remaining part of the paper, we assume that $0 < \alpha \le 1$.
Consider the following $q-$fractional initial value problem
\begin{eqnarray}\label{xds1}\left\{
\begin{array}{ll}
_{q}C^ \alpha _a y(t)= f(t,y(t)),\; a \in
\mathbb{T}_q,\\
y(a)=y_{0}.
\end{array}\right.
\end{eqnarray}
Applying $_{q}\nabla^{-\alpha}_{a} $ to both sides of (\ref{xds1}),
we obtain
\begin{equation}\label{fer}
    y(t)=y_{0}+\;_{q}\nabla^{-\alpha}_{a}f(t,y(t)).
\end{equation}
Set
\begin{equation}\label{moon}
  f(t,y(t))=x(t)y(t),
\end{equation}
where
\begin{equation}\label{SART}
  0 \le x(t) \le
\frac{1}{t^{\alpha}(1-q)^{\alpha}}.
\end{equation}
In the following, we present a comparison result for the fractional summation operator.
\begin{thm}\label{Thm1}
Let $w$ and $v$ satisfy
\begin{equation}\label{w1}
    w(t)\ge w(a)+\;_{q}\nabla^{-\alpha}_{a}x(t)w(t)
\end{equation}
and
\begin{equation}\label{v1}
     v(t)\le v(a)+\;_{q}\nabla^{-\alpha}_{a}x(t)v(t)
\end{equation}
respectively. If $w(a) \ge v(a)$, then $w(t) \ge v(t)$ for $t
\in \Lambda_{a}=\{a=q^{n_{0}},q^{n_{0}-1},\ldots\}$.
\end{thm}
\emph{Proof.} Set $u(t)=v(t)-w(t)$. We claim that $u(t) \le 0$ for
$t \in \Lambda_{a}$. Let us assume that $u(s) \le 0$ is valid for
$s=q^{n_{0}},q^{n_{0}-1},\ldots,q^{n-1}$, where $n < n_0$. Then, for $t=q^n$ we have

\begin{equation*}
u(t)=v(t)-w(t) \le [v(a)-w(a)]+
\;_{q}\nabla^{-\alpha}_{a}x(t)[v(t)-w(t)]
\end{equation*}
or
\begin{equation*}
v(t)-w(t) \le [v(a)-w(a)]+
\frac{1}{\Gamma_{q}(\alpha)}\int_{a}^{t}(t-qs)^{\alpha-1}_{q}x(s)\big(v(s)-w(s)
\big)\nabla_{q}s.
\end{equation*}
It follows that
\begin{eqnarray}\label{sou}\nonumber
    v(t)-w(t) &\le& [v(a)-w(a)]+
\frac{1}{\Gamma_{q}(\alpha)}\int_{a}^{qt}(t-qs)^{\alpha-1}_{q}x(s)\big(v(s)-w(s)
\big)\nabla_{q}s\\
&+&\frac{1}{\Gamma_{q}(\alpha)}\int_{qt}^{t}(t-qs)^{\alpha-1}_{q}x(s)\big(v(s)-w(s)
\big)\nabla_{q}s.
\end{eqnarray}
Since $v(t)-w(t) \le 0$ and $\int_{\rho(t)}^{t}f(s)\nabla
s=(t-\rho(t))f(t)$, (\ref{sou}) can be written in the form
\begin{eqnarray}\label{dfv}\nonumber
v(t)-w(t) &\le&
\frac{1}{\Gamma_{q}(\alpha)}(t-qt)(t-qt)_{q}^{\alpha-1}x(t)
\big(v(t)-w(t) \big)\\&=&(1-q)^{\alpha} t^{\alpha} x(t)
\big(v(t)-w(t) \big),
\end{eqnarray}
where
$\Gamma_{q}(\alpha)=\frac{(1-q)^{\alpha-1}_{q}}{(1-q)^{\alpha-1}}$ is used.
It follows that
$$\big(1-x(t)(1-q)^{\alpha}t^{\alpha} \big)\big(v(t)-w(t) \big) \le 0.$$
By (\ref{SART}), we conclude that $v(t)-w(t)
\le 0$.\\
\\
Define the following operator
$$_{q}\Omega_{x}\phi=~_{q}\nabla^{-\alpha}_{a}x(t)\phi(t).$$
The following lemmas are essential in the proof of the main
theorem. We only state these statements as their proofs are
straightforward.
\begin{lem}For any constant $\lambda$, we have
$$\Big\vert\;   _{q}\Omega_{\lambda}1 \Big\vert \le\; _{q}\Omega_{\vert \lambda \vert}1.$$
\end{lem}
\begin{lem}For any constant $\lambda$, we have
$$_{q}\Omega^{n}_{\lambda}1=\frac{\lambda^{n}(t-a)_q^{n \alpha}}{\Gamma_{q}(n\alpha+1)},\;\mbox{where}\;
n \in \mathbb{N}.$$
\end{lem}
\begin{lem}
Let $\lambda >0$ be such that $\vert y(t) \vert \le \lambda$ for
$t \in \Lambda_a$. Then
$$\Big\vert\; _{q}\Omega^{n}_{y} 1 \Big\vert \le\; _{q}\Omega^{n}_{\lambda} 1,\; n \in
\mathbb{N}.$$
\end{lem}
The next
result together with Theorem \ref{Thm1} will give us the desired \emph{q--}fractional Gronwall--type inequality.
\begin{thm}\label{Thm2} Let $\vert x(t) \vert \leq \frac{1}{(1-q)^\alpha t^\alpha}$ for $t \in
\Lambda_a \bigcap
[a,b]$. Then, the $q-$fractional integral equation
\begin{equation}\label{we}
    y(t)=y(a)+\;_{q}\nabla^{-\alpha}_{a}x(t)y(t)
\end{equation}
for $t \in \Lambda_a \bigcap [a,b]$ where $b \in \mathbb{R}$, has
a solution
\begin{equation}\label{gf}
    y(t)=y(a) \sum_{k=0}^{\infty}~_{q}\Omega_{x}^{k}1.
\end{equation}
\end{thm}
\emph{Proof.}  The proof is achieved by utilizing the successive
approximation method. Set
$$y_0(t)=y(a),$$
and
$$y_n(t)=y(a)+\;_{q}\nabla^{-\alpha}_{a}
x(t)y_{n-1}(t),\;\mbox{for}\;n \ge 1.$$
We observe that
$$y_1(t)=y(a)+\;_{q}\nabla^{-\alpha}_{a}
x(t)y_{0}(t)=y(a)+\;_{q}\Omega_{x}~ y(a)$$ and
$$y_2(t)=y(a)+\;_{q}\Omega_{x} \big(y(a)+\;_{q}\Omega_{x}
y(a)\big)=y(a)+\;_{q}\Omega_{x} y(a)+\;_{q}\Omega^{2}_{x}y(a).$$
Inductively, we end up with
$$y_n(t)=y(a)\sum_{k=0}^{n}\;_{q}\Omega^{k}_{x} 1,\;\;n \ge 0.$$
Taking the limit as $n \to \infty$, we have
\begin{equation}\label{xxdd}
y(t)=y(a)\sum_{k=0}^{\infty}\;_{q}\Omega^{k}_{x} 1.
\end{equation}
It remains to prove the convergence of the series in (\ref{xxdd}).
The subsequent analysis are carried out for $a=0$.

In virtue of (\ref{SART}), we obtain
\begin{eqnarray}\nonumber\label{we}
    \sum_{k=0}^{\infty}\;_{q}\Omega^{k}_{x} 1 \le
\sum_{k=0}^{\infty}\;_{q}\Omega^{k}_{\frac{1}{t^{\alpha}(1-q)^{\alpha}}}
    1 &\le& \sum_{k=0}^{\infty} \Big( \;_{q}\nabla^{-\alpha}_{0}
    \Big)^{k}\Big(\frac{1}{t^{\alpha}(1-q)^{\alpha}}\Big)\\ &\le&
    \frac{1}{(1-q)^{\alpha}}\sum_{k=0}^{\infty}\Big(\; _{q}\nabla^{-\alpha}_{0}
    \Big)^{k}(t^{-\alpha}).
\end{eqnarray}
However, for $k=1$ we observe that
$$_{q}\nabla^{-\alpha}_{0}t^{-\alpha}=\frac{\Gamma(1-\alpha)}{\Gamma(0+1)}t^{0}=\Gamma_{q}(1-\alpha).$$
For $k=2$, it follows that
$$_{q}\nabla^{-\alpha}_{0}\big(\Gamma_{q}(1-\alpha)\big)=\Gamma_{q}(1-\alpha)_{q}\nabla^{-\alpha}_{0}t^{0}=\frac{\Gamma_{q}(1-\alpha)}{\Gamma_{q}(\alpha+1)}t^{\alpha}.$$
For $k=3$, we have
$$_{q}\nabla^{-\alpha}_{0}\Big(\frac{\Gamma_{q}(1-\alpha)}{\Gamma_{q}(\alpha+1)}t^{\alpha}\Big)=\frac{\Gamma_{q}(1-\alpha)}{\Gamma_q(\alpha+1)}\frac{\Gamma_{q}(\alpha+1)}{\Gamma_{q}(\alpha+\alpha+1)}t^{2\alpha}$$
or
$$_{q}\nabla^{-\alpha}_{0}\Big(\frac{\Gamma_{q}(1-\alpha)}{\Gamma_{q}(\alpha+1)}t^{\alpha}\Big)=\frac{\Gamma_{q}(1-\alpha)}{\Gamma_q(2\alpha+1)}t^{2\alpha}.$$
For $k=4$, we get
$$_{q}\nabla^{-\alpha}_{0}\Big(\frac{\Gamma_{q}(1-\alpha)}{\Gamma_q(2\alpha+1)}t^{2\alpha}
\Big)=\frac{\Gamma_{q}(1-\alpha)}{\Gamma_{q}(3\alpha+1)}t^{3\alpha}.$$
Therefore, (\ref{we}) becomes
$$\sum_{k=0}^{\infty}\;_{q}\Omega^{k}_{x} 1 \leq
\frac{1}{(1-q)^{\alpha}}\Big[1+\Gamma_{q}(1-\alpha)+\Gamma_{q}(1-\alpha)\sum_{k=1}^{\infty}\frac{t^{k\alpha}}{\Gamma_{q}\big(k\alpha+1\big)}
      \Big].$$
Let $a_{k}=\frac{t^{(k-1)\alpha}}{\Gamma_{q}((k-1)\alpha+1)}$.
Then
\begin{eqnarray}\nonumber
\frac{a_{k}}{a_{k-1}}&=&\frac{t^{k\alpha}}{\Gamma_{q}(k\alpha+1)}\frac{\Gamma_{q}\big((k-1)\alpha+1\big)}{t^{(k-1)\alpha}}\\
&=&t^{\alpha}\frac{\Gamma_{q}\big((k-1)\alpha+1\big)}{\Gamma_{q}(k\alpha+1)}\le
\frac{\Gamma_{q}\big((k-1)\alpha+1\big)}{\Gamma_{q}(k\alpha+1)}.
\end{eqnarray}
We observe that
\begin{eqnarray}\nonumber
\frac{\Gamma_{q}\big((k-1)\alpha+1\big)}{\Gamma_{q}(k\alpha+1)}=\frac{(1-q)_{q}^{(k-1)\alpha}(1-q)^{(1-k)\alpha}}{(1-q)_{q}^{k\alpha}(1-q)^{-k\alpha}}=\frac{(1-q)_{q}^{(k-1)\alpha}(1-q)^{\alpha}}{(1-q)_{q}^{k\alpha}}.
\end{eqnarray}
Setting
$$\frac{(1-q)_{q}^{(k-1)\alpha}(1-q)^{\alpha}}{(1-q)_{q}^{k\alpha}}:=(1-q)^{\alpha}\frac{\prod_{i=0}^{\infty}\Big(\frac{1-q^{i+1}}{1-q^{i}q^{k\alpha-\alpha+1}}\Big)}{\Pi_{i=0}^{\infty}\Big(\frac{1-q^{i+1}}{1-q^{i}q^{k\alpha+1}}\Big)}$$
we deduce that
$$\lim_{k \to
\infty}(1-q)^{\alpha}\frac{\prod_{i=0}^{\infty}\Big(\frac{1-q^{i+1}}{1-q^{i}q^{k\alpha-\alpha+1}}\Big)}{\Pi_{i=0}^{\infty}\Big(\frac{1-q^{i+1}}{1-q^{i}q^{k\alpha+1}}\Big)}=(1-q)^{\alpha}\frac{\prod_{i=0}^{\infty}(1-q^{i+1})}{\Pi_{i=0}^{\infty}(1-q^{i+1})}=(1-q)^{\alpha}<1.$$
Hence, convergence is guaranteed. In case $a >0$, we can proceed in a similar way taking into account that
 $\frac{(t-a)_q^{k\alpha}}{(t-a)_q^{k\alpha-\alpha}}=(t-q^{k\alpha}q^{\alpha}a)_q^\alpha \rightarrow t^\alpha$ as $k \rightarrow \infty$.

\begin{thm} [\emph{q}--Fractional Gronwall's Lemma]\label{main}
Let $v$ and $\mu$ be nonnegative real valued functions such that
$0 \le \mu(t)<\frac{1}{t^\alpha (1-q)^\alpha}$ for all $t \in \Lambda_a$ \emph{(}in particular if $0 \le \mu(t)<\frac{1}{ (1-q)^\alpha}$\emph{)} and
$$v(t)\le v(a)+\;_{q}\nabla^{-\alpha}_{a}v(t)\mu(t).$$
Then
$$v(t)\le v(a)\sum_{k=0}^{\infty}\Omega_{\mu}^{k}1.$$
\end{thm}
The proof of the above statement is a straightforward
implementation of Theorem \ref{Thm1} and Theorem \ref{Thm2} by
setting $w(t)=v(a)\sum_{k=0}^{\infty}(\Omega_{\mu}^{k}1)(t)$.\\
\\
In case $\alpha =1$, we deduce the following  immediate
consequence of Theorem \ref{main} which can be considered as the
well known \emph{q--}Gronwall's Lemma; consult for instance the
paper \cite{wei}.
\begin{cor}
Let $0 \le \delta(t) <\frac{1}{(1-q)}$ for all $t \in \Lambda_a$.
 If
 $$v(t)\le v(a)+\int_a^{t}\delta(s) v(s)\nabla_q s.$$
 Then
 $$v(t)\le v(a)e_{q}(t,a),$$
 where $e_{q}(t,a)=~_{q}\Omega_1(1,t-a)$ is the nabla $q-$exponential function on the
 time scale $\mathbb{T}_q$.
\end{cor}

\section{Applications}
In this section we show, by the help of the $q-$fractional
Gronwall inequality proved in the previous section,  that small
changes in the initial conditions of Caputo $q-$fractional initial
value problems lead to small changes in the solution.\\
\\
Let $f(t,y)$ satisfy a Lipschitz condition with constant $0  \le
L<1$ for all $t$ and $y$.
\begin{exam}   Consider the following $q-$fractional initial value problems
\begin{eqnarray}\label{xds}\left\{
\begin{array}{ll}
_{q}\nabla^ \alpha _a \varphi(t)= f(t,\varphi(t)),\; 0 < \alpha
\le 1,\; a \in
\mathbb{T}_q,\;t \in \Lambda_{a},\\
\varphi(a)=\gamma,
\end{array}\right.
\end{eqnarray}
and
\begin{eqnarray}\label{xds}\left\{
\begin{array}{ll}
_{q}\nabla^ \alpha _a \psi(t)= f(t,\psi(t)),\; 0 < \alpha \le 1,\;
a \in
\mathbb{T}_q,\;t \in \Lambda_{a},\\
\psi(a)=\beta.
\end{array}\right.
\end{eqnarray}
It follows that
\begin{equation*}
    \varphi(t)-\psi(t)=(\gamma-\beta)+_{q}\nabla^
{-\alpha}_{a}\big[f(t,\varphi(t))-f(t,\psi(t))\big].
\end{equation*}
Taking the absolute value, we obtain
\begin{equation*}
    \vert \varphi(t)-\psi(t) \vert\leq \vert \gamma-\beta
\vert+\Big \vert ~_{q}\nabla^ {-\alpha}_{a}
f(t,\varphi(t))-f(t,\psi(t)) \Big \vert.
\end{equation*}
or
\begin{equation*}
    \vert \varphi(t)-\psi(t) \vert\leq \vert \gamma-\beta
\vert+L  ~_{q}\nabla^ {-\alpha}_{a}\vert
\varphi(t)-\psi(t)\vert.
\end{equation*}
By using Theorem \ref{main}, we get
\begin{equation*}
    \vert \varphi(t)-\psi(t) \vert\leq \vert \gamma-\beta
\vert \sum_{i=0}^{\infty}~_{q}\Omega_{L}^{i}1=\vert \gamma-\beta
\vert ~_{q}\Omega_\alpha(L,t-a).
\end{equation*}
Consider the following $q-$fractional initial value problem
\begin{eqnarray}\label{xxds}\left\{
\begin{array}{ll}
_{q}\nabla^ \alpha _a \phi(t)= f(t,\phi(t)),\; 0 < \alpha \le 1,\;
a \in
\mathbb{T}_q,\;t \in \Lambda_a\\
\phi(a)=\gamma_n,
\end{array}\right.
\end{eqnarray}
where $\gamma_n\to \gamma$. If the solution of (\ref{xxds}) is
denoted by $\phi_{n}$, then for all $t \in \Lambda_a$ we have
\begin{equation*}
    \vert \varphi(t)-\phi_n(t) \vert\leq \vert \gamma-\gamma_n
\vert \sum_{i=0}^{\infty}~_{q}\Omega_{L}^{i}1=\vert
\gamma-\gamma_n \vert ~_{q}\Omega_\alpha(L,t-a).
\end{equation*}
Hence $\vert \varphi(t)-\phi_n(t) \vert\rightarrow 0$ as
$\gamma_n\rightarrow \gamma$. This clearly verifies the dependence of solutions on the initial conditions.
\end{exam}

\end{document}